\newif\ifpreprint
\preprinttrue 

\ifpreprint
	\documentclass[a4paper, oneside]{article}
\fi

\input{preamble}

\newcommand{\TheAuthorADM}{Alberto De~Marchi}
\newcommand{\TheEmailADM}{alberto.demarchi@unibw.de}
\newcommand{\TheOrcidADM}{0000-0002-3545-6898}
\newcommand{\TheAffiliationADM}{%
	University of the Bundeswehr Munich,
	Department of Aerospace Engineering,
	Institute of Applied Mathematics and Scientific Computing,
	85577 Neubiberg/Munich,
	Germany%
}
\newcommand{\TheTitle}{A condensing approach for linear-quadratic optimization with geometric constraints}
\newcommand{\TheKeywords}{%
	Linear-quadratic optimization\KeywordsAnd
	Augmented Lagrangian\KeywordsAnd
	Geometric constraints\KeywordsAnd
	Projection oracle%
}
\newcommand{\TheAbstract}{Optimization problems with convex quadratic cost and polyhedral constraints are ubiquitous in signal processing, automatic control and decision-making. We consider here an enlarged problem class that allows to encode logical conditions and cardinality constraints, among others. In particular, we cover also situations where parts of the constraints are nonconvex and possibly complicated, but it is practical to compute projections onto this nonconvex set. Our approach combines the augmented Lagrangian framework with a solver-agnostic structure-exploiting subproblem reformulation. While convergence guarantees follow from the former, the proposed condensing technique leads to significant improvements in computational performance.}

\ifpreprint
	\title{\bfseries \TheTitle}
	\author{\TheAuthorADM\thanks{%
			\TheAffiliationADM,
			\textsc{email} \emailLink{\TheEmailADM},
			\textsc{orcid} \orcidLink{\TheOrcidADM}.%
		}%
	}%
	\date{}
\fi

\begin{document}

\maketitle

\begin{abstract}
	\TheAbstract
	
	\ifpreprint
		\medskip
		
	\fi
\end{abstract}

\section{Introduction}
This work presents a numerical approach for addressing finite-dimensional optimization problems of the form
\begin{align}\label{eq:P}
	\minimize_{x}\quad{}&
	f(x)\coloneqq\frac{1}{2} \innprod{x}{Q x} + \innprod{q}{x} \\
	\stt\quad{}& Ax \in \CC \nonumber
\end{align}
where $x\in\R^n$ are the decision variables, $f\colon \R^n\to\R$ is a convex quadratic function (defined by a symmetric matrix $Q\in\R^{n\times n}$ and a vector $q\in\R^n$), $A\in\R^{m\times n}$ is a matrix, and $\CC\subset \R^m$ a nonempty closed set.
Model \eqref{eq:P} covers classical convex quadratic programs (QP) with linear equality and inequality constraints, but it includes much more, as it allows $\CC$ be nonconvex.
In view of control applications, the class \eqref{eq:P} includes classical instances of discrete-time optimal control problems which arise, for instance, by time discretization of linear-quadratic problems of the form
\begin{align*}
	\minimize_{x(\cdot),y(\cdot),u(\cdot)}\quad{}&
	\int_{0}^{T} \|y(t)-y_{\rm ref}(t)\|^2 + \|u(t)-u_{\rm ref}(t)\|^2
	\\
	\stt\quad{}&
	\dot{x}(t) = F(x(t),u(t))
	,\quad
	y(t) = G(x(t))
	\\ &
	x(t) \in X
	,\quad
	u(t) \in U
	,\quad
	x(0) = x_{\rm init} .
\end{align*}
When the dynamics $F$ and the observation model $G$ are linear, one recovers an instance of \eqref{eq:P} upon time discretization.
Our interest in addressing \eqref{eq:P} for general, not necessarily convex, $\CC$ lies in the possibility to specify complicated constraints (nonconvex but projection-friendly) in a natural form.
Examples of sets that admit an easy-to-evaluate, possibly set-valued projection are those arising in models known as mathematical programs with equilibrium constraints (MPEC); they often result from the union of finitely many simple sets \cite{jia2023augmented,cafieri2023mixed,hall2025lcqpow}.
Notable instances are complementarity (CC), switching (SC), vanishing (VC), and either-or constraints (EOC), defined respectively by
\begin{align*}
	\text{CC}\colon\quad{}&\{ (a,b)\in\R_+^2 ~|~ ab=0 \}, \\
	\text{SC}\colon\quad{}&\{ (a,b)\in\R^2 ~|~ ab=0 \}, \\
	\text{VC}\colon\quad{}&\{ (a,b)\in\R^2 ~|~ a\geq0 ~\wedge~ ab\geq 0 \}, \\
	\text{EOC}\colon\quad{}&\{ (a,b)\in\R^2 ~|~ a\leq0 ~\vee~ b\geq 0 \}.
\end{align*}
Among others, complementarity sets appear when modelling nonsmooth dynamics \cite{stewart2010optimal} and contact mechanics \cite{menager2025contact}, while vanishing constraints play a role in mixed-integer optimal control \cite{palagachev2015mathematical,robuschi2021multiphase}.
Nonconvex sets $\CC$ naturally appear when encoding logical constraints \cite{cafieri2023mixed}, as EOC is equivalent to $a>0 \Longrightarrow b\geq 0$.
Even binary sets intersected with linear constraints can be deemed projection-friendly, since the associated projection operator boils down to solving a binary linear program,
and they can be adopted to specify combinatorial constraints and sparsity-promoting regularization costs \cite{nikitina2025hybrid}.

Unfortunately, the numerical treatment of optimization problems with disjunctive constraints is often ineffective with standard nonlinear programming techniques, since constraint qualifications are likely to fail at interesting feasible points.
Regularization \cite{hoheisel2012mathematical}, relaxation \cite{kanzow2019relaxation} and penalty \cite{hall2025lcqpow} schemes have been widely investigated and adopted to overcome this difficulty.
However, their implementation remains limited to specific constraints,
as the associated tools (such as regularization terms and penalty functions) must be tailored to the constraint set $\CC$.
Conversely,
Jia et al. \cite{jia2023augmented} study a generic approach to treat possibly complicated constraints, accessing $\CC$ only via its projection operator.
In the generality of their (nonlinear nonconvex) setting, however, it is difficult to take advantage of other structures that the problem may have.

\paragraph*{Contribution}
In this paper we build upon the generic approach of \cite{jia2023augmented} to treat the geometric constraints,
inheriting its rigorous theoretical guarantees,
while leveraging the linear-quadratic structure in \eqref{eq:P} for computational efficiency.
Our \emph{condensing} step reduces problem size and ill-conditioning.
The proposed numerical scheme poses minimal assumptions (mere convexity of $f$), is robust to redundant constraints, and handles the constraint set $\CC$ only through its projection oracle.
The advantages (faster convergence, scalability and greater flexibility) are showcased across a range of control tasks.
Further speed-up is demonstrated when exploiting also ``safe'' equality constraints, e.g., those arising from dynamics.

\paragraph*{Outline}
The notion of solution is discussed in \cref{sec:optimality} before describing the augmented Lagrangian framework in \cref{sec:methods}.
The condensing approach and related procedure are detailed in \cref{sec:minimize_over_x}.
Numerical results and comparisons are reported in \cref{sec:numerics}.

\paragraph*{Notation}
Given a set $S\subseteq\R^m$,
its \emphdef{indicator} is defined by $\indicator_S(v)=0$ if $v\in S$, and $\indicator_S(v)=\infty$ otherwise.
The associated \emphdef{projection}
$\proj_S\colon \R^m\rightrightarrows \R^m$ is defined by $\proj_S(v) \coloneqq \argmin_{z\in S}\|z-v\|$,
while $\dist_S(v) \coloneqq \inf_{z\in S}\|z-v\|$ indicates the \emphdef{distance} of $v$ from $S$.
The identity matrix $\identity$ has suitable size, clear from the context.

\section{Optimality and Stationarity}\label{sec:optimality}

Problem \eqref{eq:P} can be rewritten in the equivalent form
\begin{equation}\label{eq:Pz}
	\minimize_{x ,\, z}~
	f(x) + \indicator_\CC(z)
	\qquad
	\stt~
	Ax - z = 0
\end{equation}
by introducing the auxiliary variable $z\in\R^m$.
Then, the \emphdef{Lagrangian function} for \eqref{eq:Pz} reads
\begin{equation*}
	\LL(x,z,y)
	\coloneqq
	f(x)
	+ \indicator_\CC(z)
	+ \innprod{y}{Ax-z}
\end{equation*}
where $y\in\R^m$ denotes the Lagrange multiplier associated to the equality constraints.
By differentiation one obtains the following stationarity conditions, necessary for first-order optimality (under suitable constraint qualifications):
\begin{subequations}\label{eq:KKT_Pz}
	\begin{align}
		&0 = \nabla_x \LL(x,z,y) = \nabla f(x) + A^\top y , \\
		&0 \in \partial_z \LL(x,z,y) = \limnormalcone_\CC(z) - y , \\
		&0 = \nabla_y \LL(x,z,y) = Ax - z ,
	\end{align}
\end{subequations}
where $\partial$ and $\limnormalcone_\CC$ denote the limiting subdifferential and the limiting normal cone to $\CC$,
respectively, as in \cite{jia2023augmented}.
The ``limiting'' character appears as $\CC$ is potentially nonconvex.
Formally eliminating the auxiliary $z$, the stationarity conditions for the original problem \eqref{eq:P} read
\begin{align*}
	0 ={}& \nabla f(x) + A^\top y , &
	y \in{}& \limnormalcone_\CC(Ax) ,
\end{align*}
where the inclusion encodes both feasibility of $x$, as it implies $Ax\in\CC$, and complementarity of $y$ with the constraint.\footnote{%
	Observe that, for any $\alpha>0$,
	the condition $Ax\in\proj_{\CC}(Ax + \alpha y)$ implies $y \in \limnormalcone_\CC(Ax)$.
	In contrast, due to nonconvexity of $\CC$, the converse implication
	$y \in \limnormalcone_\CC(Ax)
	\Longrightarrow
	Ax\in\proj_{\CC}(Ax+\alpha y)$
	is valid only for sufficiently small $\alpha>0$.
	For convex $\CC$, they are equivalent for all $\alpha>0$.
}
A point $x \in\R^n$ is called a \emphdef{Fritz-John} (FJ) point for \eqref{eq:P} if there exists $(\mu, y) \in\R\times\R^m$ with $(\mu, y) \neq 0$ such that
\begin{align}\label{eq:FJ_P}
	0 ={}& \mu \nabla f(x) + A^\top y , &
	y \in{}& \limnormalcone_\CC(Ax) .
\end{align}
A FJ point is a \emphdef{Karush-Kuhn-Tucker} (KKT) point whenever $\mu \neq 0$.
In that case, $y/\mu$ is the associated vector of Lagrange multipliers.
A FJ point for which $\mu = 0$ is called a \emphdef{singular stationary} point:
a feasible point at which constraint qualifications fail (e.g., LICQ does not hold).
A point $x\in\R^n$ is called an \emphdef{infeasible stationary} point of problem \eqref{eq:P} if $Ax \not\in \CC$ and $A^\top \left[ Ax-\proj_{\CC}(Ax) \right] = 0$:
an infeasible point for \eqref{eq:P} that is a stationary point for the feasibility problem $\minimize_{x} \dist_{\CC}^2\left( Ax \right)$.

For implementing an iterative numerical procedure, a termination condition can be based on an approximate notion of optimality.
Given some tolerances $\epsilon_{\rm d},\epsilon_{\rm p}>0$, a point $x\in\R^n$ is called \emphdef{$(\epsilon_{\rm d},\epsilon_{\rm p})$-stationary} (in the FJ sense)
if there exists $(\mu,y,z)\in\R\times\R^m\times\CC$ with $(\mu,y)\neq 0$ such that the conditions
\begin{subequations}
	\label{eq:approxKKT_Pz}
	\begin{align}
		\| \mu \nabla f(x) + A^\top y \| \leq{}& \epsilon_{\rm d} \\
		\dist(y, \limnormalcone_\CC(z)) \leq{}& \epsilon_{\rm d} \\
		\| Ax - z \|_\infty \leq{}& \epsilon_{\rm p}
	\end{align}
\end{subequations}
are satisfied, see \eqref{eq:KKT_Pz} and \eqref{eq:FJ_P}.
Tolerances $\epsilon_{\rm d}$ and $\epsilon_{\rm p}$ reflect the relaxation on \emph{d}ual and \emph{p}rimal residuals, respectively.

\section{Augmented Lagrangian Method}\label{sec:methods}

The auxiliary variable $z\in\R^m$ has been introduced in \eqref{eq:Pz} to derive the stationarity conditions for \eqref{eq:P} and to define a suitable termination criterion based on approximate optimality.
However, the expression of problem \eqref{eq:P} in the splitting form \eqref{eq:Pz} is also closer to the numerical procedure described below, as it allows us to favorably decouple the linear operator $A$ from the set $\CC$.
We will relax and penalize violation of the linear constraint $Ax=z$, following the augmented Lagrangian (AL) framework \cite{jia2023augmented,birgin2014practical}.
This entails solving a sequence of subproblems,
at least inexactly and up to stationarity (in contrast to global optimality),
and then updating some parameters based on the progress thus obtained.
A key role is played by the augmented Lagrangian function for \eqref{eq:Pz}, here defined for all $x\in\R^n$, $z\in\CC$ by
\begin{equation}
	\label{eq:augLagr}
	\LL_\sigma(x,z,\widehat{x},\widehat{y})
	\coloneqq
	\paramf f(x)
	+ \frac{1}{2}\|Ax-z+\widehat{y}\|^2
	+ \frac{\paramx}{2} \|x-\widehat{x}\|^2
\end{equation}
for a given primal-dual estimate $(\widehat{x},\widehat{y})\in\R^n\times\R^m$ and a pair of positive penalty parameters $\sigma \coloneqq (\paramx,\paramf)$.
Note that, for reasons of numerical stability, we prefer to scale down the cost function $f$ instead of increasing the quadratic constraint penalization, following \cite{armand2019augmented}, while also adding a proximal term for Tikhonov-like primal regularization.

\subsection{Safeguarded augmented Lagrangian scheme}

The AL approach requires solving subproblems of the form
\begin{align}
	\label{eq:ALMsubproblem_xz}
	\tag{AL$_k$}
	\minimize_{x ,\, z}~{}&
	\LL_{\sigma_k}(x,z,\widehat{x}_k,\widehat{y}_k) &
	\stt~{}&
	z\in\CC ,
\end{align}
up to approximate stationarity,
to obtain the next primal estimate $(x_k,z_k)$.
The overall scheme is outlined in \cref{alg:ALM}, where it is clear that the minimization of \eqref{eq:ALMsubproblem_xz} is the most computationally demanding task.
To alleviate this effort, warm-starting \cref{step:ALM:x} from the current estimate $(\widehat{x}_k,\widehat{z}_k)$ is a common practice.
After the primal update, \cref{step:ALM:y} generates a dual estimate
$
	y_k \gets \widehat{y}_k + A x_k - z_k
$
by matching the subproblem's stationarity conditions with the KKT conditions for \eqref{eq:Pz}.
Note that, at \cref{step:ALM:hat}, we adopt the safeguarding scheme of \cite[\S 4]{birgin2014practical} to guarantee global convergence:
the previous dual estimate is projected onto a compact (arbitrarily large) set $Y\subset \R^m$ to obtain
$\widehat{y}^k$.
Then, \cref{step:ALM:termination} tests if the current primal-dual point is $(\epsilon_{\rm d},\epsilon_{\rm p})$-stationary, returning when successful.
Owing to the $\varepsilon_k$-stationarity conditions for \eqref{eq:ALMsubproblem_xz} and the dual update at \cref{step:ALM:y}, every triplet $(x_k,z_k,y_k)$ satisfies
\begin{align*}
	\left\| \paramf_{k} \nabla f(x_k) + A^\top y_k + \paramx_{k} (x_k-\widehat{x}_k) \right\|
	\leq{}&
	\varepsilon_k \\
	\dist\left( y_k, \limnormalcone_\CC(z_k) \right)
	\leq{}&
	\varepsilon_k
	,
\end{align*}
meaning that it suffices to check the computable quantities $E_k$ and $V_k$ defined at \cref{step:ALM:E,step:ALM:V} to assert approximate optimality of $x_k$, in the sense of \eqref{eq:approxKKT_Pz}.

Finally, the penalty parameters and tolerance controlling the inner minimization are updated.
If a sufficient reduction of the infeasibility measure $V_k$ has been obtained,
both penalty parameters are kept constant and the inner tolerance is decreased at \cref{step:update_subtol}.
Otherwise, the latter remains constant and the penalty parameters are reduced at \cref{step:update_penalty}.

\begin{algorithm2e}[tbh]
	\DontPrintSemicolon%
	\KwIn{$x^0\in\R^n$, $y^0\in\R^m$, $z^0\in\R^m$, $\epsilon_{\rm p},\epsilon_{\rm d}>0$}
	\KwData{$\kappa_V\in[0,1)$, $\kappa_\varepsilon, \kappa_{\paramf} \in(0,1)$, $\kappa_{\paramx}\in(0,1]$, $\varepsilon_1,\paramf_{1},\paramx_{1}>0$, $Y\subset\R^m$ compact}
	$V_0\gets \infty$\;
	\For{$k = 1,2\ldots$}{
		$(\widehat{x}_k, \widehat{z}_k, \widehat{y}_k) \gets \left(x_{k-1}, z_{k-1}, \proj_Y(y_{k-1})\right)$\label{step:ALM:hat}\;
		$(x_k,z_k) \gets \varepsilon_k$-stationary point of \eqref{eq:ALMsubproblem_xz}\label{step:ALM:x}\;
		$y_k \gets \widehat{y}_k + A x_k-z_k$\label{step:ALM:y}\;
		$E_k \gets \max\{ \|\mu_k \nabla f(x_k) + A^\top y_k\|, \varepsilon_k \}$\label{step:ALM:E}\;
		$V_k \gets \|A x_k-z_k\|_\infty$\label{step:ALM:V}\;
		\lIf{\normalfont$E_k\leq \epsilon_{\rm d}$ and $V_k\leq \epsilon_{\rm p}$\label{step:ALM:termination}}{%
			\KwRet{$(x_k,z_k,y_k)$}
		}
		\If{$V_k \leq \max\{\epsilon_{\rm p},\kappa_V V_{k-1}\}$}{%
			$\paramx_{k+1}, \paramf_{k+1}, \varepsilon_{k+1} \gets \paramx_{k}, \paramf_{k}, \kappa_\varepsilon\max\{\epsilon_{\rm d}, \varepsilon_k\}$\label{step:update_subtol}\;
		}
		\Else{%
			$\paramx_{k+1}, \paramf_{k+1}, \varepsilon_{k+1} \gets \kappa_{\paramx} \paramx_{k}, \kappa_{\paramf} \paramf_{k}, \varepsilon_k$\label{step:update_penalty}\;
		}%
	}
	\caption{Augmented Lagrangian scheme}
	\label{alg:ALM}
\end{algorithm2e}

Before discussing how the subproblem at \cref{step:ALM:x} can be solved with off-the-shelf methods,
let us review some convergence guarantees available for \cref{alg:ALM}.

\subsection{Convergence and infeasibility detection}

An essential property of \cref{alg:ALM} is that subproblem \eqref{eq:ALMsubproblem_xz} is well-posed for every $\sigma_k>0$, $\widehat{x}_k\in\R^n$ and $\widehat{y}_k\in\R^m$.
This fact follows from the strict convexity of the augmented Lagrangian $\LL_{\sigma_k}(\cdot,\cdot,\widehat{x}_k,\widehat{y}_k)$ and closedness of set $\CC$.
Then, classical convergence guarantees are directly inherited.
The main observation is that feasible limit points are automatically KKT points (under constraints qualifications).
\begin{mybox}
	\begin{proposition}
		Suppose the sequence $\{\varepsilon_k\}$ in \cref{alg:ALM} satisfies $\varepsilon_k\to 0$.
		Then each feasible accumulation point $(\bar{x},\bar{z})$ of a sequence $\{(x_k,z_k)\}$ generated by \cref{alg:ALM} is an asymptotically KKT point of \eqref{eq:Pz}.
		Moreover, if matrix $A$ has full rank, each feasible accumulation point $\bar{x}$ of $\{x_k\}$ is a KKT point of \eqref{eq:P}.
	\end{proposition}
\end{mybox}
\begin{proof}
	See \cite[Thm 4.3]{jia2023augmented}, \cite[Thm 4.1]{armand2019augmented}.
\end{proof}
Since the feasible set \eqref{eq:P} is nonconvex, it may happen that a limit point is not feasible,
but even in this case it contains some useful information:
infeasible limit points are stationary for the constraint violation.
\begin{mybox}
	\begin{proposition}
		Suppose the sequence $\{\varepsilon_k\}$ in \cref{alg:ALM} is bounded.
		Then each accumulation point $\bar{x}$ of a sequence $\{x_k\}$ generated by \cref{alg:ALM} is a stationary point for the so-called feasibility problem
		\begin{equation*}
			\minimize_{x}\quad \dist_{\CC}^2(Ax) .
		\end{equation*}
	\end{proposition}
\end{mybox}
\begin{proof}
	See \cite[Prop. 4.2]{jia2023augmented}, \cite[Thm 6.3]{birgin2014practical}.
\end{proof}

\subsection{Inner solver: projected gradient}

The subproblems \eqref{eq:ALMsubproblem_xz} arising at \cref{step:ALM:x} of \cref{alg:ALM} have the general form
\begin{equation}\label{eq:generic_subproblem}
	\minimize_w~
	\varphi(w)
	\qquad
	\stt~
	w\in W
\end{equation}
with a continuously differentiable function $\func{\varphi}{\mathbb{W}}{\R}$ and some nonempty, closed set $W\subset\mathbb{W}$, where $\mathbb{W}$ is an arbitrary Euclidean space.
The projection operator associated to $W$ is not required to be single-valued, but it should be easy to evaluate in practice.
One can use
$w\coloneqq (x,z)$, $W \coloneqq \R^n\times \CC$, and function $\varphi\coloneqq w\mapsto \LL_{\sigma_k}(x,z,\widehat{x}_k,\widehat{y}_k)$
to recover \eqref{eq:ALMsubproblem_xz}.

The prototypical method for addressing \eqref{eq:generic_subproblem} is known as \emph{projected gradient},
for which neither $\varphi$ nor $W$ need to be convex \cite[\S 3]{jia2023augmented}.
Iterates are constructed based on the update rule
\begin{equation*}
	w_{j+1}
	\gets
	\proj_W \left( w_j - \gamma_j \nabla \varphi(w_j) \right) ,
\end{equation*}
where $\gamma_j>0$ is a suitable stepsize, typically obtained by backtracking linesearch to deliver sufficient decrease in the cost $\varphi$.
The projection operation ensures that each iterate is feasible, namely $w_j\in W$ for all $j\in\N$. 
Projected, or semismooth, Newton-type methods can achieve fast convergence but are applicable only to simple sets $W$ in practice (such as hyperboxes).
In contrast, first-order schemes such as projected gradient are more flexible and can benefit from nonmonotone globalization mechanisms \cite{jia2023augmented,demarchi2023monotony} as well as (quasi-Newton) extrapolation techniques \cite{demarchi2022proximal}.

\section{Condensing: How and Why}\label{sec:minimize_over_x}

Although subproblem \eqref{eq:ALMsubproblem_xz} can be solved with off-the-shelf solvers without additional effort,
its special structure is particularly favorable and should be exploited.
Since \cref{step:ALM:x} has the highest computational cost in \cref{alg:ALM},
the overall solver performance can improve significantly.
We demonstrate how to take advantage of it, on top of warm-starting the inexact solves at \cref{step:ALM:x}, before assessing its beneficial effects with numerical examples.

Our starting point here is the observation that, for fixed $z\in\CC$,
\eqref{eq:ALMsubproblem_xz} becomes an unconstrained problem in $x$ with a strictly convex quadratic objective:
\begin{equation}\label{eq:marginal_x_subproblem}
	\minimize_{x}~{}
	\paramf f(x)
	+ \frac{1}{2}\|Ax-z+\widehat{y}\|^2
	+ \frac{\paramx}{2} \|x-\widehat{x}\|^2
	.
\end{equation}
Therefore, for any given $z\in\CC$ there is a unique optimal $x$ to minimize $\LL_{\sigma_k}(\cdot,z,\widehat{x}_k,\widehat{y}_k)$;
we denote $\func{\XX_k}{\CC}{\R^n}$ this single-valued mapping.
How to efficiently evaluate $\XX_k$ will be the topic of \cref{sec:compute_marginal_function}.
With this in mind, we can also define the (marginal) function $\func{\VV_k}{\CC}{\R}$ as
\begin{align}
	\VV_k(z)
	\coloneqq{}&
	\min_x
	\LL_{\sigma_k}(x,z,\widehat{x}_k,\widehat{y}_k)
	\nonumber\\
	={}&
	\LL_{\sigma_k}(\XX_k(z),z,\widehat{x}_k,\widehat{y}_k)
	\label{eq:marginal_cost_z}
	,
\end{align}
which gives the equivalent, yet \emph{condensed} AL subproblem
\begin{equation}
	\label{eq:ALMsubproblem_z}
	\tag{CAL$_k$}
	\minimize_{z}~{}
	\VV_k(z)
	\qquad
	\stt~{}
	z\in\CC .
\end{equation}
Inspired by the ``implicit'' approach in \cite{zheng2019unified},
the formal minimization over the original variable $x$ leaves us with a subproblem in the auxiliary variable $z$ only.
Note that \eqref{eq:ALMsubproblem_z} is again of the form \eqref{eq:generic_subproblem},
and it is smaller in size than \eqref{eq:ALMsubproblem_xz} (with $m$ variables instead of $n+m$),
as well as better conditioned \cite[Thm 1]{zheng2019unified}.

Crucially in practice, $\VV_k$ is continuously differentiable and has, in fact, even globally Lipschitz continuous gradient.%
\arxivonly{Proof in \cref{app:condensing_globally_lipschitz}.}
Moreover, its value and gradient are relatively cheap to compute, as we illustrate below.
\eqref{eq:ALMsubproblem_z} is the \emph{condensed} counterpart of \eqref{eq:ALMsubproblem_xz} and can be addressed with the same off-the-shelf tools.
Therefore, one can easily replace \cref{step:ALM:x} of \cref{alg:ALM} with
\begin{itemize}
	\item $z_k \gets \varepsilon_k$-stationary point of \eqref{eq:ALMsubproblem_z},
	\item $x_k\gets \XX_k(z_k)$.
\end{itemize}

\subsection{Evaluating the marginal function}\label{sec:compute_marginal_function}

By construction, $\XX_k$ returns the unique minimizer $x$ of $\LL_{\sigma_k}(\cdot,z,\widehat{x}_k,\widehat{y}_k)$, given $z\in\CC$.
Owing to convexity, this corresponds to solving $0 = \nabla_x \LL_{\sigma_k}(x,z,\widehat{x}_k,\widehat{y}_k)$,
which can be rewritten as the linear system
\begin{equation}\label{eq:linsys_dense}
	\begin{bmatrix}
		\paramf_k Q + \paramx_k \identity + A^\top A
	\end{bmatrix}
	x
	= \begin{pmatrix}
		\paramx_k \widehat{x}_k - \paramf_k q + A^\top (z-\widehat{y}_k)
	\end{pmatrix}
	.
\end{equation}
Note that the matrix on the left-hand side is always symmetric positive definite, hence nonsingular,
and independent on $z$, whereas the right-hand side vector varies linearly in $z$.
Thus, $\XX_k$ is a linear map from $\CC$ to $\R^n$.
Moreover, factorization caching can be adopted,
effectively factorizing (at most) once for each solve of \eqref{eq:ALMsubproblem_z}.
To avoid the matrix fill-in due to $A^\top A$,
in our implementation we do not solve \eqref{eq:linsys_dense} as is, but lift it to recover structure;
see \cref{sec:hard_equality_constraints} below for details.

Given $z$, one can find the optimal $\XX_k(z)$ and then evaluate $\VV_k(z)$ directly through its definition \eqref{eq:marginal_cost_z}.
Owing to the chain rule, the gradient $\nabla \VV_k(z)$ is readily obtained as
\begin{equation*}
	\nabla \VV_k(z)
	={}
	\nabla_z \LL_{\sigma_k}(\XX_k(z),z,\widehat{x}_k,\widehat{y}_k)
	={}
	z - A \XX_k(z) - \widehat{y}_k
	,
\end{equation*}
since $\nabla_x \LL_{\sigma_k}(\XX_k(z),z,\widehat{x}_k,\widehat{y}_k) = 0$ by definition of $\XX_k$.

\subsection{Lower-level equality constraints}\label{sec:hard_equality_constraints}

The formulation \eqref{eq:P} offers a compact yet generic model to analyse.
However, practical optimization models often exhibit a richer structure:
linear dynamics in optimal control tasks give rise to linear equality constraints, say $\Aeq x = \beq$,
that are ``safe'', in the sense that they are feasible and non-redundant.
Here we proceed to demonstrate how this additional structure can be exploited for speeding up \cref{step:ALM:x} even further.
Convergence guarantees remain valid, as demonstrated in \cite[Chapter 6]{birgin2014practical}.

By keeping the equality constraints $\Aeq x = \beq$ without relaxation,
we treat them as lower-level constraints in the AL subproblem (on par with the inclusion $z\in\CC$)
instead of including them in the augmented Lagrangian.
Then, the counterpart of \eqref{eq:marginal_x_subproblem} is an equality-constrained quadratic program,
\begin{equation*}
	\XX_k(z)
	\coloneqq{}
	\arg\min_x
	\left\{
	\LL_{\sigma_k}(x,z,\widehat{x}_k,\widehat{y}_k)
	\,\middle\vert\,
	\Aeq x = \beq
	\right\},
\end{equation*}
whose solution requires only a linear system.
Upon lifting with auxiliary variables $\lambda$ and $\lambda_{\rm eq}$,
we write it as
\begin{equation*}
	\begin{bmatrix}
		\paramf_k Q + \paramx_k \identity & A^\top & \Aeq^\top \\
		A & - \identity & 0 \\
		\Aeq & 0 & 0
	\end{bmatrix}
	\begin{pmatrix}
		x \\ \lambda \\ \lambda_{\rm eq}
	\end{pmatrix} = \begin{pmatrix}
		\paramx_k \widehat{x}_k - \paramf_k q \\
		z - \widehat{y}_k \\
		\beq
	\end{pmatrix}
\end{equation*}
to retain sparsity and symmetry.%
\arxivonly{See \cref{app:condensing_linsys} for details.}

\section{Numerical Results}\label{sec:numerics}

We consider three different benchmarks:
an initial value problem of a discontinuous dynamic with single switch,
a discretized obstacle problem,
and a reference tracking task with the AFTI-16 aircraft model.

We assess various solver configurations along two axes:
\begin{itemize}
	\item subproblem:
	\emph{extended} formulation (as in \cite{jia2023augmented}) against the proposed \emph{condensed} one (using the marginalization of \cref{sec:compute_marginal_function}).
	In the latter case, explicit equality constraints $\Aeq x =\beq$ are treated either as \emph{soft} (relaxing them like $Ax=z$) or \emph{hard} (according to \cref{sec:hard_equality_constraints}).
	\item subsolver:
	\nmpg{} (with spectral stepsize) \cite{demarchi2023monotony} or \panocp{} (with quasi-Newton extrapolation) \cite{demarchi2022proximal}, both with nonmonotone linesearch.
\end{itemize}
Overall, six solver configurations are compared, with subproblem formulations that increasingly exploit the problem's structure.

\paragraph*{Solver setup}
Our prototypical implementation in MATLAB 2025a run on a standard laptop, Intel Core i7 2.80 GHz, 16 GB RAM.
For each call to \cref{alg:ALM}, we set a time limit of 100 s,
tolerances $\epsilon_{\rm d}=\epsilon_{\rm p}=10^{-6}$,
and parameters
$\kappa_{\paramx} = 1$,
$\paramx_1 = 10^{-6}$,
$\kappa_{\paramf} = \nicefrac{1}{4}$,
$\paramf_1 = 1$,
$\kappa_\varepsilon = \nicefrac{1}{2}$,
$\varepsilon_1 = 1$,
$\kappa_V = \nicefrac{9}{10}$,
$Y \coloneqq [-10^{20},10^{20}]^m$.
For each problem instance, we consider 10 random (normally distributed) initializations for $x^0$, and then set $y^0 = 0$, $z^0 = \proj_{\CC}(A x^0)$.

\paragraph*{Profiles}
Solvers are compared by means of data and (extended) performance profiles.
For $P$ the set of problem instances and $S$ the set of solver settings, let $t_{s,p}$ denote the wall-clock runtime required by $s \in S$ on $p \in P$.
When solver $s$ fails on problem $p$, we set $t_{s,p}=\infty$.
A \emph{data profile} depicts the empirical distribution $f_s^{\rm abs} \colon [0, \infty) \to [0, 1]$ of $t_{s,\cdot}$ over $P$, namely
$f_s^{\rm abs}(t) \coloneqq |\{p \in P ~|~ t_{s,p} \leq t\}| / |P|$,
where $|P|$ denotes the cardinality of $P$.
As such, each data profile reports the fraction of problems $f_s^{\rm abs}(t)$ that can be solved (for given tolerances) by solver $s$ with a computational budget $t$,
and therefore it is independent of the other solvers.
In contrast, an \emph{extended performance profile} is relative to a pool of solvers $S$;
it depicts the empirical distribution $f_s^{\rm rel} \colon [0, \infty) \to [0, 1]$ of the ratio $r_{s,\cdot}$ of $t_{s,\cdot}$ with the best runtime among the other solvers, namely
$r_{s,p} \coloneqq t_{s,p} / \min \{ t_{s^\prime,p} \colon s^\prime\in S ,\, s^\prime\ne s \}$.%
\arxivonly{Extended performance profiles were introduced by Mahajan, Leyffer \& Kirches in \S 4 of their technical report ``Solving mixed-integer nonlinear programs by QP-diving'' (2012), ANL/MCS Preprint P2071-0312, Argonne National Laboratory.}

\subsection{Initial value problem}
Introduced by Stewart and Anitescu \cite[\S 2]{stewart2010optimal},
this test example contains a dynamical system with a discontinuous right-hand side, in which a single switch occurs.
The discontinuous dynamics
can be reformulated into a linear complementarity system.
We then discretize the system using the implicit Euler scheme with $N$ nodes, as in \cite[\S 5.2]{hall2025lcqpow}, obtaining
\begin{align}
	\label{eq:initialValueProblem}
	\minimize_{\substack{x_0,\ldots,x_N\in\R\\ y_1,\ldots,y_N\in\R\\ \lambda_1,\ldots,\lambda_N\in\R}}\quad{}& (x_N-x_{\rm ref})^2 + h \sum_{k=0}^{N-1} x_k^2 \\
	\stt\quad{}& x_k = x_{k-1} + h (3-2y_k) ,\nonumber\\
	& 0 \leq x_k + \lambda_k \perp 1-y_k \geq 0 ,\nonumber\\
	& 0 \leq \lambda_k \perp y_k \geq 0 \qquad\qquad \text{for}~k=1,\ldots,N \nonumber
\end{align}
with reference $x_{\rm ref} \coloneqq \nicefrac{5}{3}$ and
time step $h = 2/N$.

\begin{figure}[tbh]
	\centering%
	\includegraphics{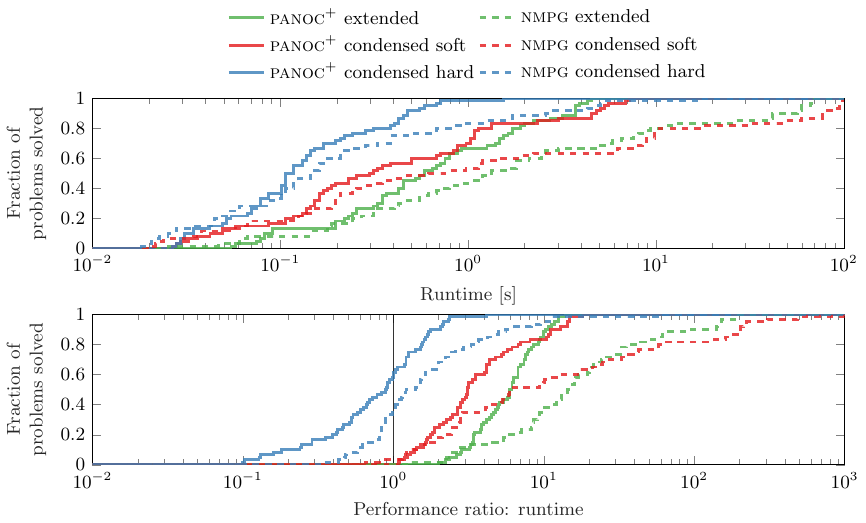}%
	\caption{Initial value problem \eqref{eq:initialValueProblem}: data (top) and performance (bottom) profiles for comparing different solver settings.}%
	\label{fig:ivp}%
\end{figure}

Results with increasing number of discretization intervals $N \in \{2^3,2^4,\ldots,2^8\}$, for a total of 60 calls to each solver, are reported in \cref{fig:ivp}.
All profiles indicate that the \emph{condensed} formulation beats the extended,
with the \emph{hard} treatment of equality constraints delivering better runtimes.
For fixed subproblem, \panocp{} typically outperforms \nmpg{} as subsolver.
Moreover, we point out that, despite the generic approach powering \cref{alg:ALM}, it achieves performance on par with the solver
developed in \cite{hall2025lcqpow} specifically for complementarity constraints.

\subsection{Obstacle problem}

We consider the optimal control of a discretized obstacle problem as investigated in
\cite[Example 6.2]{jia2023augmented}.
The problem is given by
\begin{align}
	\label{eq:obstacleProblem}
	\minimize_{x,y,z\in\R^N}\quad{}& \frac{1}{2}\|x\|^2 + \frac{1}{2}\|y\|^2 - \innprod{1}{y} \\
	\stt\quad{}& x \geq 0 ,\quad
	\min\{y,z\} = 0 ,\quad
	x+Ay=z \nonumber
\end{align}
where $A\in\R^{N\times N}$ is a tridiagonal matrix that arises from a discretization of the negative Laplace operator in one dimension.
Its simple solution $x=y=z=0$ is degenerate and makes the problem particularly arduous to solve.

\begin{figure}[tbh]
	\centering%
	\includegraphics{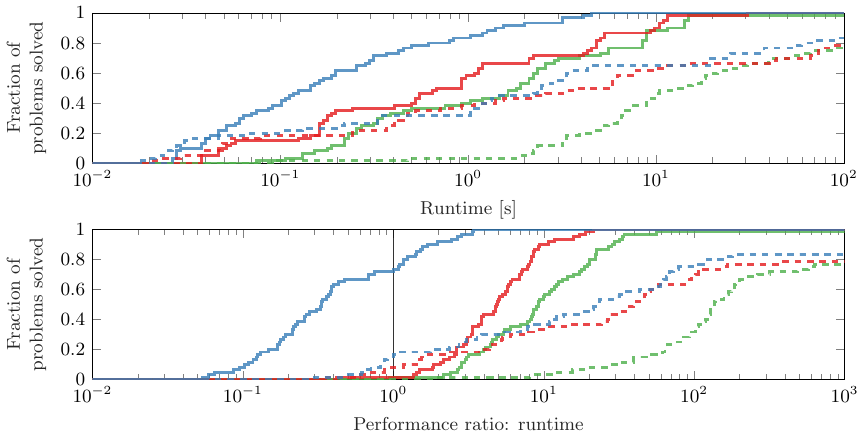}%
	\caption{Obstacle problem \eqref{eq:obstacleProblem}: data (top) and performance (bottom) profiles for comparing different solver settings. Legend as in \cref{fig:ivp}.}%
	\label{fig:obstacle}%
\end{figure}

The numerical results depicted in \cref{fig:obstacle} lead to observations analogous to those for \cref{fig:ivp}.
Compared to \cite[\S 6.1]{jia2023augmented}, we can push the discretization from their $N=64$ to $N=256$ and, using \panocp{} as subsolver, keep a reasonable runtime and the total number of inner iterations in the order of thousands instead of millions.

\subsection{AFTI-16 tracking control}\label{sec:afti16}

The AFTI-16 model from \cite{kapasouris1990design}
provides a LTI approximation of the aircraft's longitudinal dynamics at 3000 ft altitude and 0.6 Mach speed:
the state description includes forward velocity, angle of attack, pitch rate, and pitch angle,
while control inputs are elevator and flaperon angles.
The model, converted to discrete time by sampling
with zero-order hold on the inputs and a sample time of $T_{\rm s} = 50$ ms,%
\arxivonly{More details in \cref{app:afti16}.}
is $x_{k+1} = A_{\rm DT} x_k + B_{\rm DT} u_k$ with
\begin{equation*}
	A_{\rm DT}
	=
	\begin{bmatrix}
		0.9993 & -3.0083 & -0.1131 & -1.6081 \\
		-4.703 \cdot 10^{-6} & 0.9862 & 0.0478 & 3.85 \cdot 10^{-6} \\
		3.703 \cdot 10^{-6} & 2.0833 & 1.0089 & -4.362 \cdot 10^{-6} \\
		1.356 \cdot 10^{-7} & 0.0526 & 0.0498 & 1
	\end{bmatrix}
	,\,
	B_{\rm DT}
	=
	\begin{bmatrix}
		-0.08045 & -0.6347 \\
		-0.02914 & -0.01428 \\
		-0.8679 & -0.0913 \\
		-0.02159 & -0.002181
	\end{bmatrix}
	,\,
	C_{\rm DT}^\top
	=
	\begin{bmatrix}
		0 & 0 \\
		1 & 0 \\
		0 & 0 \\
		0 & 1
	\end{bmatrix}
	.
\end{equation*}
The control task is to track the null reference for the observed states, with an additional cost for the control effort:
\begin{align}\label{eq:afti16problem}
	\minimize_{\substack{x_1,\ldots,x_N\in\R\\ u_0,\ldots,u_{N-1}\in\R}}\quad{}& \sum_{k=0}^{N-1} \| C_{\rm DT} x_{k+1} \|^2 + \left\| \frac{u_k}{\bar{u}} \right\|^2 \\
	\stt\quad{}& x_{k+1} = A_{\rm DT} x_k + B_{\rm DT} u_k , \nonumber\\
	& u_k \in U \qquad\qquad \text{for}~k=0,\ldots,N-1 . \nonumber
\end{align}
The initial state $x_0 = x_{\rm init}$ is fixed, as well as the time horizon $N$.
As detailed in \eqref{eq:afti16problem}, the control system is also subject to a geometric constraint on the admissible controls.
Set $U \coloneqq \{(a,b)\in\R^2 \colon |a|\leq \bar{u}, |b|\leq \bar{u}, ab=0\}$
encodes an either-or condition intersected with simple bounds:
at any given time, at most one control input can be applied, and its magnitude cannot be larger than $\bar{u} \coloneqq 25$.

\begin{figure}[tbh]
	\centering%
	\includegraphics{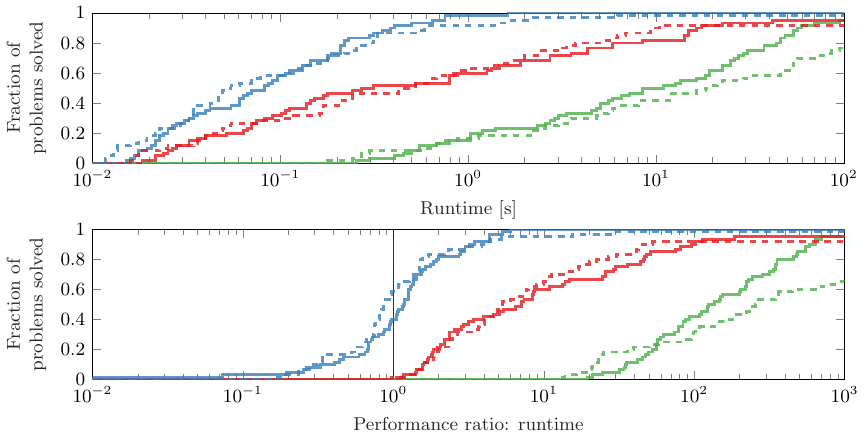}%
	\caption{AFTI-16 problem \eqref{eq:initialValueProblem}: data (top) and performance (bottom) profiles for comparing different solver settings. Legend as in \cref{fig:ivp}.}%
	\label{fig:afti16}%
\end{figure}

The outcomes depicted in \cref{fig:afti16} (obtained with randomized $x_{\rm init}$) demonstrate an even greater advantage for the condensed formulation with hard equality constraints.
Moreover, compared to previous examples, the performance in this case seems less dependent on the subsolver.

\begin{figure}[tbh]
	\centering%
	\includegraphics{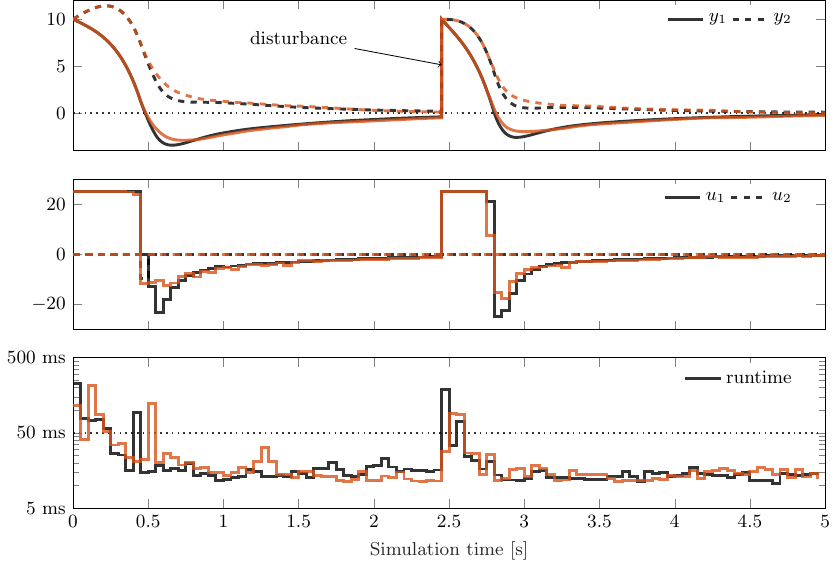}%
	\caption{Closed-loop simulation of AFTI-16 behavior: outputs (top), inputs (middle) and solver runtime (bottom) using cold- (black) or warm-starting (orange) between consecutive simulation steps. A disturbance is injected at half time, resetting the observed states as at the initial time.}%
	\label{fig:afti16_mpc_plant}%
\end{figure}

To assess the quality of solutions,
we
repeatedly solve \eqref{eq:afti16problem} in a simulated environment with closed-loop control;
here we start from $x_0=(10,10,10,10)$ and use the solver configuration ``\nmpg{} condensed hard''.
Numerical results with and without warm-starting the solver are reported in \cref{fig:afti16_mpc_plant}.
Despite the nonconvex admissible set $U$, the solver quickly finds sufficiently good control sequences to steer the system as desired,
also counteracting the disturbance injected at $t=2.5$.

\section{Conclusions}

In this article, we introduced a structure-exploiting approach for addressing linear-quadratic optimization problems with geometric constraints.
Our strategy inherits convergence guarantees from the augmented Lagrangian framework and builds upon simple and efficient computational kernels.
Making analytical progress in the subproblem formulation, we effectively precondition it and reduce its size to advantage of the user-defined subsolver.
Extensive numerical experiments demonstrate significant speed-up and improved robustness.%
\arxivonly{A scalability analysis is provided in \cref{app:more_num_plots}.}

\ifpreprint
	\phantomsection
	\addcontentsline{toc}{section}{References}%
	\bibliographystyle{habbrv}
	\bibliography{biblio}
\fi

\ifpreprint
\appendix
\section{Additional Material}

\subsection{Globally Lipschitz gradient}\label{app:condensing_globally_lipschitz}

Let us denote $\|\XX_k\|$ the operator norm of $\XX_k$,
which depends on the problem data as well as on parameters $\sigma_k$.
Now we proceed to show that, for each $k\in\N$, $\VV_k$ has globally Lipschitz continuous gradient, regardless of the problem data and choice of parameters $\sigma_k$.
Picking two arbitrary $z_1,z_2\in\R^m$ and invoking the triangle and Cauchy-Schwartz inequalities, we have the bound
\begin{align*}
	\| \nabla \VV_k(z_1) - \nabla \VV_k(z_2) \|
	&=
	\| \nabla_z \LL_{\sigma_k}(\XX_k(z_1),z_1,\widehat{x}_k,\widehat{y}_k) - \nabla_z \LL_{\sigma_k}(\XX_k(z_2),z_2,\widehat{x}_k,\widehat{y}_k) \|
	\\
	&=
	\| (z_1 - A \XX_k(z_1) - \widehat{y}_k) - (z_2 - A \XX_k(z_2) - \widehat{y}_k) \|
	\\
	&=
	\| A [\XX_k(z_2) - \XX_k(z_1)] + (z_1 - z_2) \|
	\\
	&\leq
	\| A [\XX_k(z_2) - \XX_k(z_1)] \| + \| z_1 - z_2 \|
	\\
	&\leq
	\|A\| \| \XX_k(z_2) - \XX_k(z_1) \| + \| z_1 - z_2 \|
	\\
	&\leq
	(\|A\| \|\XX_k\| + 1) \| z_1 - z_2 \| .
\end{align*}
Since $\XX_k$ is a linear operator (also in the presence of equalities $\Aeq x = \beq$), the finite quantity $L_{\VV_k}\coloneqq\|A\| \|\XX_k\| + 1$ represents a global Lipschitz constant for $\nabla \VV_k$.

\subsection{Structured linear system}\label{app:condensing_linsys}

Let us derive step-by-step the linear system in \cref{sec:hard_equality_constraints}.
Given some $z\in\CC$, the optimal $x=\XX_k(z)$ is obtained as the unique solution to
\begin{align*}
	\minimize_x\quad{}&
	\paramf_k f(x) + \frac{1}{2}\|Ax-z+\widehat{y}_k\|^2 + \frac{\paramx_k}{2} \|x-\widehat{x}_k\|^2 \\
	\stt\quad{}&
	\Aeq x = \beq. \nonumber
\end{align*}
The optimality conditions for this problem read
\begin{align*}
	0 ={}& \paramf_k (Qx+q) + A^\top (Ax-z+\widehat{y}_k) +\Aeq^\top \lambda_{\rm eq}  + \paramx_k(x-\widehat{x}_k)\\
	0 ={}& \Aeq x - \beq
\end{align*}
where $\lambda_{\rm eq}$ denotes the Lagrange multiplier for $\Aeq x = \beq$.
Introducing the auxiliary variable $\lambda \equiv Ax - z + \widehat{y}_k$, we can write the equivalent system of linear equations
\begin{align*}
	0 ={}& \paramf_k (Qx+q) + A^\top \lambda + \Aeq^\top \lambda_{\rm eq} + \paramx_k(x-\widehat{x}_k) \\
	0 ={}& Ax - z + \widehat{y}_k - \lambda \\
	0 ={}& \Aeq x - \beq ,
\end{align*}
which can be rearranged in matrix form obtaining the one in \cref{sec:hard_equality_constraints}.

\subsection{AFTI-16 model}\label{app:afti16}

The continuous-time (dynamics and observation) model from \cite{kapasouris1990design} is $\dot{x} = A_{\rm CT} x + B_{\rm CT} u$, $y = C_{\rm CT} x$ with
\begin{equation*}
	A_{\rm CT} ={} \begin{bmatrix}
			-0.0151 & -60.5651 & 0 & -32.174 \\
			-0.0001 & -1.3411 & 0.9929 & 0 \\
			0.00018 & 43.2541 & -0.86939 & 0 \\
			0 & 0 & 1 & 0
		\end{bmatrix}
		,\,
	B_{\rm CT} ={} \begin{bmatrix}
			-2.516 & -13.136 \\
			-0.1689 & -0.2514 \\
			-17.251 & -1.5766 \\
			0 & 0
		\end{bmatrix}
		,\,
	C_{\rm CT}^\top
		=
		\begin{bmatrix}
		0 & 0 \\
		1 & 0 \\
		0 & 0 \\
		0 & 1
		\end{bmatrix}
		.
\end{equation*}
The discrete-time counterpart given in \cref{sec:afti16} was obtained with MATLAB's Control System Toolbox via the commands
\begin{itemize}
	\item [\texttt{>>}] \texttt{plantCT = ss(A\_CT, B\_CT, C\_CT, D\_CT)};
	\item [\texttt{>>}] \texttt{plantDT = c2d(plantCT, Ts)};
\end{itemize} 
and taking the state-space matrix representation of \texttt{plantDT} with four decimals for each entry.
Both models have $D_{\rm CT} = D_{\rm DT} = 0$ and $C_{\rm CT} = C_{\rm DT} \ne 0$.

\subsection{Scalability outcomes}\label{app:more_num_plots}

In addition to the data and performance profiles presented in \cref{sec:numerics},
we report here \emph{scalability profiles}.
These display the median values and the interquartile range over all instances of a given size, indicating how the metric is affected by the problem size.

Scalability profiles are depicted in \cref{fig:N_vs_runtime}.
It appears that, among all solver configurations,
the condensed, hard formulation with \panocp{} subsolver consistently performs best and scales better with problem size.
The solver settings with \panocp{} are significantly more effective than with \nmpg{} for problems \eqref{eq:initialValueProblem} and \eqref{eq:obstacleProblem}, but they are surprisingly similar for \eqref{eq:afti16problem}.

\begin{figure}[tbh]
	\centering%
	\includegraphics{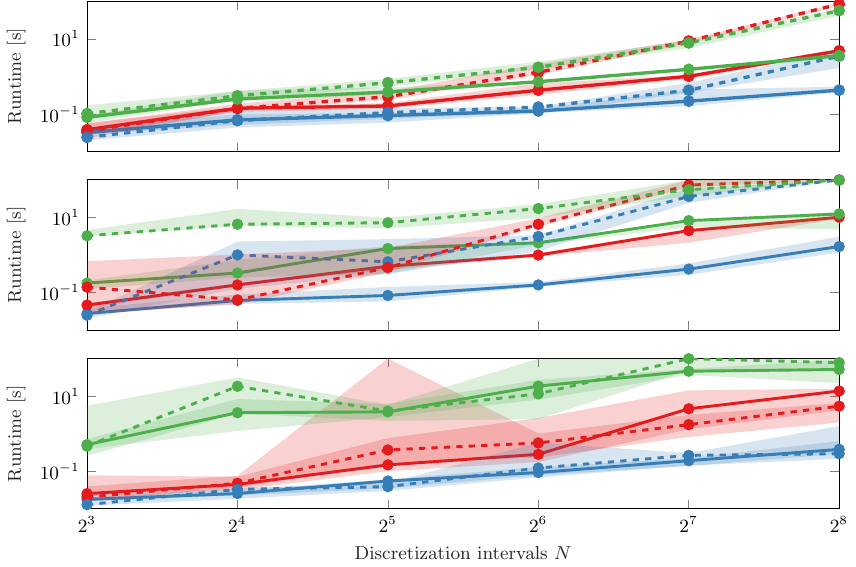}%
	\caption{Scalability profiles for comparing different solver settings: initial value problem \eqref{eq:initialValueProblem} (top), obstacle problem \eqref{eq:obstacleProblem} (middle) and AFTI-16 problem \eqref{eq:afti16problem} (bottom). Legend as in \cref{fig:ivp}.}%
	\label{fig:N_vs_runtime}%
\end{figure}
\fi

\end{document}